\documentclass{amsart}
\usepackage{amstext}
\usepackage{amsthm}
\usepackage{amssymb}

\makeatletter
\numberwithin{equation}{section}
\numberwithin{figure}{section}
\theoremstyle{plain}
\newtheorem{thm}{Theorem}
  \theoremstyle{plain}
  \newtheorem{cor}[thm]{Corollary}

\allowdisplaybreaks
\usepackage[sort&compress,square,numbers]{natbib}

\makeatother

\begin{document}

\title{Symmetric identities for an analogue of Catalan polynomials}

\author{Taekyun Kim}
\address{Department of Mathematics, Kwangwoon University, Seoul 139-701, Republic
of Korea}
\email{tkkim@kw.ac.kr}

\author{Dae San Kim}
\address{Department of Mathematics, Sogang University, Seoul 121-742, Republic
of Korea}
\email{dskim@sogang.ac.kr}

\author{Jong-Jin Seo}
\address{Department of Applied Mathematics, Pukyong National University, Busan, Republic
of Korea}
\email{seo2011@pknu.ac.kr}

\begin{abstract}
In this paper, we consider an analogue of Catalan polynomials and
give some identities of symmetry for those polynomials by using fermionic $p$-adic integrals on the ring of $p$-adic integers.
\end{abstract}

\maketitle
\global\long\def\relphantom#1{\mathrel{\phantom{{#1}}}}

\section{Introduction}

Let $p$ be an odd prime number. Throughout this paper, $\mathbb{Z}_{p}$,
$\mathbb{Q}_{p}$ and $\mathbb{C}_{p}$ will denote the ring of $p$-adic
integers, the field of $p$-adic
rational numbers, and the completion of the algebraic closure of $\mathbb{Q}_{p}$.
The $p$-adic norm $\left|\cdot\right|_{p}$ is normalized as $\left|p\right|_{p}=\frac{1}{p}$.
Let $f$ be a continuous function on $\mathbb{Z}_{p}$.
Then the fermionic $p$-adic integral of $f$ on $\mathbb{Z}_{p}$ is defined
by Kim
\begin{align}
\int_{\mathbb{Z}_{p}}f\left(x\right)d\mu_{-1}\left(x\right) & =\lim_{N\rightarrow\infty}\sum_{x=0}^{p^{N}-1}f\left(x\right)\mu_{-1}\left(x\right)\label{eq:1}\\
 & =\lim_{N\rightarrow\infty}\sum_{x=0}^{p^{N}-1}f\left(x\right)\left(-1\right)^{x},\quad\left(\text{see \cite{6, 7, 8, 9, 10, 13, 14}}\right).\nonumber
\end{align}

From (\ref{eq:1}), we note that
\begin{equation}
\int_{\mathbb{Z}_{p}}f_{1}\left(x\right)d\mu_{-1}\left(x\right)+\int_{\mathbb{Z}_{p}}f\left(x\right)d\mu_{-1}\left(x\right)=2f\left(0\right),\label{eq:2}
\end{equation}
where $f_{1}\left(x\right)=f\left(x+1\right)$, (see \cite{6, 7, 8, 9, 10, 13, 14}).

By (\ref{eq:2}), we get
\begin{equation}
\int_{\mathbb{Z}_{p}}f\left(x+n\right)d\mu_{-1}\left(x\right)+\left(-1\right)^{n-1}\int_{\mathbb{Z}_{p}}f\left(x\right)d\mu_{-1}\left(x\right)=2\sum_{l=0}^{n-1}\left(-1\right)^{n-1-l}f\left(l\right),\label{eq:3}
\end{equation}
where $n\in\mathbb{N}$, (see \cite{6, 7, 8, 9, 10, 13, 14}).

As is well known, the Catalan numbers are defined by the generating
function to be
\begin{equation}
\frac{1-\sqrt{1-4t}}{2t}=\frac{2}{1+\sqrt{1-4t}}=\sum_{n=0}^{\infty}C_{n}t^{n}, \text{(see \cite{1, 2, 3, 4, 5, 11, 12, 15, 16}),}\label{eq:4}
\end{equation}
where $C_{n}=\binom{2n}{n}\frac{1}{n+1}$, $\left(n\ge0\right)$.

In addition, the Catalan polynomials are given by the
generating function
\begin{equation}
\frac{2}{1+\sqrt{1-4t}}\left(1-4t\right)^{\frac{x}{2}}=\sum_{n=0}^{\infty}C_{n}\left(x\right)t^{n}, \text{(see \cite{1, 2, 3, 4, 5, 11, 12, 15, 16})}.\label{eq:5}
\end{equation}

Thus, by (\ref{eq:4}) and (\ref{eq:5}), we get
\begin{equation}
C_{n}\left(x\right)=\sum_{m=0}^{n}\sum_{j=0}^{m}\left(\frac{x}{2}\right)^{j}S_{1}\left(m,j\right)\left(-4\right)^{m}\frac{C_{n-m}}{m!},\label{eq:6}
\end{equation}
where $S_{1}\left(m,n\right)$ is the Stirling number of the first
kind.

For $w\in\mathbb{N}$, we define $w$-Catalan polynomials which are
given by the generating function
\begin{equation}
\frac{2}{1+\left(1-4t\right)^{\frac{w}{2}}}\left(1-4t\right)^{\frac{w}{2}x}=\sum_{n=0}^{\infty}C_{n,w}\left(x\right)t^{n}.\label{eq:7}
\end{equation}

From (\ref{eq:7}), we note that $C_{n,1}\left(x\right)=C_{n}\left(x\right)$.

When $x=0$, $C_{n,w}=C_{n,w}\left(0\right)$ are called $w$-Catalan
numbers.

From (\ref{eq:2}), we note that
\begin{align}
\int_{\mathbb{Z}_{p}}\left(1-4t\right)^{\frac{w}{2}\left(x+y\right)}d\mu_{-1}\left(y\right) & =\frac{2}{\left(1-4t\right)^{\frac{w}{2}}+1}\left(1-4t\right)^{\frac{w}{2}x}\label{eq:8}\\
 & =\sum_{n=0}^{\infty}C_{n,w}\left(x\right)t^{n},\nonumber
\end{align}
where $t\in\mathbb{C}_{p}$ with $\left|t\right|_{p}<p^{-\frac{1}{p-1}}$.

In particular, for $w=1$, we have
\begin{align}
\int_{\mathbb{Z}_{p}}\left(1-4t\right)^{\frac{x+y}{2}}d\mu_{-1}\left(y\right) & =\frac{2}{1+\left(1-4t\right)^{\frac{1}{2}}}\left(1-4t\right)^{\frac{x}{2}}\label{eq:9}\\
 & =\sum_{n=0}^{\infty}C_{n}\left(x\right)t^{n}.\nonumber
\end{align}

For $k\ge0$, $d,w\in\mathbb{N}$ with $d\equiv1\pmod{2}$, $w\equiv1\pmod{2}$,
we define the function $S_{k,d}\left(w-1\right)$ as follows:
\begin{equation}
S_{k,d}\left(w-1\right)=\sum_{i=0}^{w-1}\binom{\frac{di}{2}}{k}\left(-1\right)^{i}.\label{eq:10}
\end{equation}
When $d=1$,
\[
S_{k,1}\left(w-1\right)=S_{k}\left(w-1\right)=\sum_{i=0}^{w-1}\binom{\frac{i}{2}}{k}\left(-1\right)^{i}.
\]

In this paper, we give some identities of symmetry for the $w$-Catalan
polynomials which are derived from the fermionic $p$-adic integrals
on $\mathbb{Z}_{p}$.

\section{Identities of symmetry for the $w$-Catalan polynomials}

From (\ref{eq:9}), we note that
\begin{equation}
\int_{\mathbb{Z}_{p}}\binom{\frac{x+y}{2}}{n}d\mu_{-1}\left(y\right)=\frac{\left(-1\right)^{n}}{4^{n}}C_{n}\left(x\right),\quad\left(n\ge0\right).\label{eq:11}
\end{equation}

For $d\in\mathbb{N}$ with $d\equiv1\pmod{2}$, we have
\begin{align}
 & \int_{\mathbb{Z}_{p}}\left(1-4t\right)^{\frac{x+d}{2}}d\mu_{-1}\left(x\right)+\int_{\mathbb{Z}_{p}}\left(1-4t\right)^{\frac{x}{2}}d\mu_{-1}\left(x\right)\label{eq:12}\\
 & =2\sum_{i=0}^{d-1}\left(1-4t\right)^{\frac{i}{2}}\left(-1\right)^{i}\nonumber \\
 & =\sum_{n=0}^{\infty}\left(2\sum_{i=0}^{d-1}\binom{\frac{i}{2}}{n}\left(-4\right)^{n}\left(-1\right)^{i}\right)t^{n}.\nonumber
\end{align}

Thus, by (\ref{eq:9}) and (\ref{eq:12}), we get
\begin{equation}
\sum_{n=0}^{\infty}\left(C_{n}\left(d\right)+C_{n}\right)t^{n}=\sum_{n=0}^{\infty}\left(2^{2n+1}\sum_{i=0}^{d-1}\binom{\frac{i}{2}}{n}\left(-1\right)^{n-i}\right)t^{n}.\label{eq:13}
\end{equation}

Therefore, by (\ref{eq:13}), we obtain the following theorem.
\begin{thm}
\label{thm:1} For $n\ge0$, $d\in\mathbb{N}$ with $d\equiv1\pmod{2}$,
we have
\[
C_{n}\left(d\right)+C_{n}=2^{2n+1}\sum_{i=0}^{d-1}\binom{\frac{i}{2}}{n}\left(-1\right)^{n-i}.
\]

\end{thm}
Now, we observe that, for any $d\in\mathbb{N}$,
\begin{align}
 & \int_{\mathbb{Z}_{p}}\left(1-4t\right)^{\frac{x+d}{2}}d\mu_{-1}\left(x\right)+\int_{\mathbb{Z}_{p}}\left(1-4t\right)^{\frac{x}{2}}d\mu_{-1}\left(x\right)\label{eq:14}\\
 & =\frac{2}{\int_{\mathbb{Z}_{p}}\left(1-4t\right)^{\frac{d}{2}x}d\mu_{-1}\left(x\right)}\int_{\mathbb{Z}_{p}}\left(1-4t\right)^{\frac{x}{2}}d\mu_{-1}\left(x\right).\nonumber
\end{align}

From (\ref{eq:8}), we have
\begin{equation}
\int_{\mathbb{Z}_{p}}\left(1-4t\right)^{\frac{d}{2}x}d\mu_{-1}\left(x\right)=\frac{2}{1+\left(1-4t\right)^{\frac{d}{2}}}=\sum_{n=0}^{\infty}C_{n,d}t^{n}. \label{eq:15}
\end{equation}

Let $w_{1},w_{2}\in\mathbb{N}$ with $w_{1}\equiv1\pmod{2}$, $w_{2}\equiv1\pmod{2}$.

By using double fermionic $p$-adic invariant integral on $\mathbb{Z}_{p}$,
we get
\begin{align}
 & \frac{\int_{\mathbb{Z}_{p}}\int_{\mathbb{Z}_{p}}\left(1-4t\right)^{\frac{w_{1}x+w_{2}y}{2}}d\mu_{-1}\left(x\right)d\mu_{-1}\left(y\right)}{\int_{\mathbb{Z}_{p}}\left(1-4t\right)^{\frac{w_{1}w_{2}}{2}x}d\mu_{-1}\left(x\right)}\label{eq:16}\\
 & =\frac{2\left(\left(1-4t\right)^{\frac{w_{1}w_{2}}{2}}+1\right)}{\left(\left(1-4t\right)^{\frac{w_{1}}{2}}+1\right)\left(\left(1-4t\right)^{\frac{w_{2}}{2}}+1\right)}.\nonumber
\end{align}

Now, we also consider the following fermionic $p$-adic invariant
integral on $\mathbb{Z}_{p}$ associated with $w$-Catalan polynomials.

\begin{align}
 & I\label{eq:17}\\
 & =\frac{\int_{\mathbb{Z}_{p}}\int_{\mathbb{Z}_{p}}\left(1-4t\right)^{\frac{w_{1}x_{1}+w_{2}x_{2}+w_{1}w_{2}x}{2}}d\mu_{-1}\left(x_{1}\right)d\mu_{-1}\left(x_{2}\right)}{\int_{\mathbb{Z}_{p}}\left(1-4t\right)^{\frac{w_{1}w_{2}}{2}x}d\mu_{-1}\left(x\right)}\nonumber \\
 & =\frac{2\left(1-4t\right)^{\frac{w_{1}w_{2}}{2}x}\left(\left(1-4t\right)^{\frac{w_{1}w_{2}}{2}}+1\right)}{\left(\left(1-4t\right)^{\frac{w_{1}}{2}}+1\right)\left(\left(1-4t\right)^{\frac{w_{2}}{2}}+1\right)}.\nonumber
\end{align}

From (\ref{eq:13}), (\ref{eq:14}), we note that, $w_{1}, d\in\mathbb{N}$ with $w_{1}\equiv1\pmod{2}$,
\begin{align}
\frac{\int_{\mathbb{Z}_{p}}\left(1-4t\right)^{\frac{dx}{2}}d\mu_{-1}\left(x\right)}{\int_{\mathbb{Z}_{p}}\left(1-4t\right)^{\frac{dw_{1}x}{2}}d\mu_{-1}\left(x\right)} & =\sum_{i=0}^{w_{1}-1}\left(1-4t\right)^{\frac{di}{2}}\left(-1\right)^{i}\label{eq:18}\\
 & =\sum_{k=0}^{\infty}\left(\sum_{i=0}^{w_{1}-1}\binom{\frac{di}{2}}{k}\left(-1\right)^{i}\left(-4\right)^{k}\right)t^{k}\nonumber \\
 & =\sum_{k=0}^{\infty}S_{k,d}\left(w_{1}-1\right)\left(-4\right)^{k}t^{k}.\nonumber
\end{align}

By (\ref{eq:17}) and (\ref{eq:18}), we get
\begin{align}
 & I\label{eq:19}\\
 & =\left(\int_{\mathbb{Z}_{p}}\left(1-4t\right)^{\frac{w_{1}\left(x_{1}+w_{2}x\right)}{2}}d\mu_{-1}\left(x_{1}\right)\right)\left(\frac{\int_{\mathbb{Z}_{p}}\left(1-4t\right)^{\frac{w_{2}x_{2}}{2}}d\mu_{-1}\left(x_{2}\right)}{\int_{\mathbb{Z}_{p}}\left(1-4t\right)^{\frac{w_{1}w_{2}}{2}x}d\mu_{-1}\left(x\right)}\right)\nonumber \\
 & =\left(\sum_{l=0}^{\infty}C_{l,w_{1}}\left(w_{2}x\right)t^{l}\right)\left(\sum_{k=0}^{\infty}S_{k,w_{2}}\left(w_{1}-1\right)\left(-4\right)^{k}t^{k}\right)\nonumber \\
 & =\sum_{n=0}^{\infty}\left(\sum_{l=0}^{n}\left(-4\right)^{n-l}C_{l,w_{1}}\left(w_{2}x\right)S_{n-l,w_{2}}\left(w_{1}-1\right)\right)t^{n}.\nonumber
\end{align}

On the other hand,
\begin{align}
 & I\label{eq:20}\\
 & =\left(\int_{\mathbb{Z}_{p}}\left(1-4t\right)^{\frac{w_{2}\left(x_{2}+w_{1}x\right)}{2}}d\mu_{-1}\left(x_{2}\right)\right)\left(\frac{\int_{\mathbb{Z}_{p}}\left(1-4t\right)^{\frac{w_{1}x_{1}}{2}}d\mu_{-1}\left(x_{1}\right)}{\int_{\mathbb{Z}_{p}}\left(1-4t\right)^{\frac{w_{1}w_{2}}{2}x}d\mu_{-1}\left(x\right)}\right)\nonumber \\
& =\left(\sum_{l=0}^{\infty}C_{l,w_{2}}\left(w_{1}x\right)t^{l}\right)\left(\sum_{k=0}^{\infty}S_{k,w_{1}}\left(w_{2}-1\right)\left(-4\right)^{k}t^{k}\right)\nonumber \\
 & =\sum_{n=0}^{\infty}\left(\sum_{l=0}^{n}\left(-4\right)^{n-l}C_{l,w_{2}}\left(w_{1}x\right)S_{n-l,w_{1}}\left(w_{2}-1\right)\right)t^{n}.\nonumber
\end{align}

Therefore, by (\ref{eq:19}) and (\ref{eq:20}), we obtain the following
theorem.
\begin{thm}
\label{thm:2} For $n\ge0$, $w_{1},w_{2}\in\mathbb{N}$ with $w_{1}\equiv1\pmod{2}$
and $w_{2}\equiv1\pmod{2}$, we have
\[
\sum_{l=0}^{n}\left(-4\right)^{n-l}C_{l,w_{1}}\left(w_{2}x\right)S_{n-l,w_{2}}\left(w_{1}-1\right)=\sum_{l=0}^{n}\left(-4\right)^{n-l}C_{l,w_{2}}\left(w_{1}x\right)S_{n-l,w_{1}}\left(w_{2}-1\right).
\]

\end{thm}
Let $w_{1}=1$. Then we have
\begin{equation}
C_{n,1}\left(w_{2}x\right)=\sum_{l=0}^{n}\left(-4\right)^{n-l}C_{l,w_{2}}\left(x\right)S_{n-l,1}\left(w_{2}-1\right).\label{eq:21}
\end{equation}

Therefore, by (\ref{eq:21}), we obtain the following corollary.
\begin{cor}
\label{cor:3}For $w_{2}\in\mathbb{N}$ with $w_{2}\equiv1\pmod{2}$,
we have
\[
C_{n}\left(w_{2}x\right)=\sum_{l=0}^{n}\left(-4\right)^{n-l}C_{l,w_{2}}\left(x\right)S_{n-l}\left(w_{2}-1\right).
\]

\end{cor}
When $x=0$ in Theorem \ref{thm:2}, we have the following corollary.
\begin{cor}
\label{cor:4} For $n\ge0$, $w_{1},w_{2}\in\mathbb{N}$ with $w_{1}\equiv1\pmod{2}$,
$w_{2}\equiv1\pmod{2}$, we have

\[
\sum_{l=0}^{n}\left(-4\right)^{n-l}C_{l,w_{1}}S_{n-l,w_{2}}\left(w_{1}-1\right)=\sum_{l=0}^{n}\left(-4\right)^{n-l}C_{l,w_{2}}S_{n-l,w_{1}}\left(w_{2}-1\right).
\]

\end{cor}

From (\ref{eq:17}), we have

\begin{align}
 & I\label{eq:22}\\
 & =\left(\left(1-4t\right)^{\frac{w_{1}w_{2}}{2}x}\int_{\mathbb{Z}_{p}}\left(1-4t\right)^{\frac{w_{1}x_{1}}{2}}d\mu_{-1}\left(x_{1}\right)\right)\left(\frac{\int_{\mathbb{Z}_{p}}\left(1-4t\right)^{\frac{w_{2}x_{2}}{2}}d\mu_{-1}\left(x_{2}\right)}{\int_{\mathbb{Z}_{p}}\left(1-4t\right)^{\frac{w_{1}w_{2}}{2}x}d\mu_{-1}\left(x\right)}\right)\nonumber \\
 & =\left(\left(1-4t\right)^{\frac{w_{1}w_{2}}{2}x}\int_{\mathbb{Z}_{p}}\left(1-4t\right)^{\frac{w_{1}x_{1}}{2}}d\mu_{-1}\left(x_{1}\right)\right)\left(\sum_{l=0}^{w_{1}-1}\left(-1\right)^{l}\left(1-4t\right)^{\frac{w_{2}}{2}l}\right)\nonumber \\
 & =\sum_{l=0}^{w_{1}-1}\left(-1\right)^{l}\int_{\mathbb{Z}_{p}}\left(1-4t\right)^{\frac{w_{1}}{2}\left(x_{1}+w_{2}x+\frac{w_{2}}{w_{1}}l\right)}d\mu_{-1}\left(x_{1}\right)\nonumber \\
 & =\sum_{n=0}^{\infty}\left(\sum_{l=0}^{w_{1}-1}\left(-1\right)^{l}C_{n,w_{1}}\left(w_{2}x+\frac{w_{2}}{w_{1}}l\right)\right)t^{n}.\nonumber
\end{align}

On the other hand,
\begin{align}
 & I\label{eq:23}\\
 & =\left(\left(1-4t\right)^{\frac{w_{1}w_{2}}{2}x}\int_{\mathbb{Z}_{p}}\left(1-4t\right)^{\frac{w_{2}x_{2}}{2}}d\mu_{-1}\left(x_{2}\right)\right)\left(\frac{\int_{\mathbb{Z}_{p}}\left(1-4t\right)^{\frac{w_{1}x_{1}}{2}}d\mu_{-1}\left(x_{1}\right)}{\int_{\mathbb{Z}_{p}}\left(1-4t\right)^{\frac{w_{1}w_{2}}{2}x}d\mu_{-1}\left(x\right)}\right)\nonumber \\
 & =\left(\left(1-4t\right)^{\frac{w_{1}w_{2}}{2}x}\int_{\mathbb{Z}_{p}}\left(1-4t\right)^{\frac{w_{2}x_{2}}{2}}d\mu_{-1}\left(x_{2}\right)\right)\left(\sum_{l=0}^{w_{2}-1}\left(-1\right)^{l}\left(1-4t\right)^{\frac{w_{1}}{2}l}\right)\nonumber \\
 & =\sum_{l=0}^{w_{2}-1}\left(-1\right)^{l}\int_{\mathbb{Z}_{p}}\left(1-4t\right)^{\frac{w_{2}}{2}\left(x_{2}+w_{1}x+\frac{w_{1}}{w_{2}}l\right)}d\mu_{-1}\left(x_{2}\right)\nonumber \\
 & =\sum_{n=0}^{\infty}\left(\sum_{l=0}^{w_{2}-1}\left(-1\right)^{l}C_{n,w_{2}}\left(w_{1}x+\frac{w_{1}}{w_{2}}l\right)\right)t^{n}.\nonumber
\end{align}

Therefore, by (\ref{eq:22}) and (\ref{eq:23}), we obtian the following
theorem.
\begin{thm}
\label{thm:5} For $n\ge0$, $w_{1},w_{2}\in\mathbb{N}$ with $w_{1}\equiv1\pmod{2}$,
$w_{2}\equiv1\pmod{2}$, we have
\[
\sum_{l=0}^{w_{1}-1}\left(-1\right)^{l}C_{n,w_{1}}\left(w_{2}x+\frac{w_{2}}{w_{1}}l\right)=\sum_{l=0}^{w_{2}-1}\left(-1\right)^{l}C_{n,w_{2}}\left(w_{1}x+\frac{w_{1}}{w_{2}}l\right).
\]

\end{thm}
Setting $w_{2}=1$ in Theorem \ref{thm:5}, we get the multiplication--
type formula for the Catalan polynomials as follows:
\[
C_{n}\left(w_{1}x\right)=\sum_{l=0}^{w_{1}-1}\left(-1\right)^{l}C_{n,w_{1}}\left(x+\frac{1}{w_{1}}l\right).
\]
\bibliographystyle{amsplain}

\begin{thebibliography}{10}
\bibitem{1} H. W. Gould,
\emph{Sums and convolved sums of Catalan numbers and their generating functions,
} Indian J. Math. 46 (2004), no. 2-3, 137-160.

\bibitem{2} R. Hampel,
\emph{On the problem of Catalan,
} (Polish) Prace Mat. 4 (1960), 11-19.

\bibitem{3} Y. He and C. Wang,
\emph{Recurrence formulae for Apostol-Bernoulli and Apostol-Euler polynomials
}Adv. Difference Equ. 2012, 2012:209, 16 pp.

\bibitem{4} S. Hyyro,
\emph{ On the Catalan problem,
}(Finnish) Arkhimedes 1963 (1963), no. 1, 53-54.

\bibitem{5} K. Inkeri,
\emph{On Catalan's problem,
}  Acta Arith. 9 (1964), 285-290.

\bibitem{6} D. S. Kim and K. H. Park, \textit{Identities of symmetry for Euler polynomials arising   from quotients of fermionic integrals invariant under
$S_{3}$}, J. Inequal. Appl. 2010, Art. ID \textbf{851521}, 16 pp.

\bibitem{7} T. Kim, \textit{Symmetry $p$-adic invariant integral on $\mathbb Z_{p}$ for Bernoulli and Euler polynomials}, J. Difference Equ. Appl. \textbf{14} (2008), 1267-1277.

\bibitem{8} T. Kim, \textit{Symmetry of power sum polynomials and multivariate fermionic $p$-adic invariant integral on $\mathbb{Z}_p$}. Russ. J. Math. Phys. 1616 (2009), no.\textbf{1}, 93-96.

\bibitem{9} T. Kim, \textit{$q$-Volkenborn integration}, Russ. J. Math. Phys. 9 (2002), no.\textbf{3}, 288-299.

\bibitem{10} T. Kim, \textit{On $p$-adic interpolating function for $q$-Euler numbers and its derivatives}, J. Math. Anal. Appl. 339 (2008), no. \textbf{1}, 598-608.

\bibitem{11} A. Natucci,
\emph{Ricerche sistematiche intorno al "teorema di Catalan",
} (Italian) Giorn. Mat. Battaglini (5) 2(82) (1954), 297-300.

\bibitem{12} H. Ozden, I. Naci Cangul and Y. Simsek,\textit{Multivariate interpolation functions of higher-order $q$-Euler numbers and their applications}, Abstract and
Applied Analysis, Article Number: \textbf{390857} Published: 2008.

\bibitem{13} H. Ozden and Y. Simsek, \textit{A new extension of $q$-Euler numbers and polynomials related to their interpolation functions}, Applied Mathematics Letters, Volume \textbf{21}, Issue 9, September 2008, Pages 934-939.

\bibitem{14} R. Rangarajan and P. Shashikala,
\emph{A pair of classical orthogonal polynomials connected to Catalan numbers,
} Adv. Stud. Contemp. Math. (Kyungshang) 23 (2013), no. 2, 323-335.

\bibitem{15} A. Rotkiewicz,
\emph{ Sur le probl\'{e}me de Catalan,
} (French) Elem. Math. 15 (1960), 121-124.

\bibitem{16} A. D. Sands,
\emph{On generalised Catalan numbers,
} Discrete Math. 21 (1978), no 2. 219-221.


\end{thebibliography}

\end{document}